\documentclass[preprints,article,accept,moreauthors,pdftex]{my-mdpi} 

\usepackage{amsmath,amsthm,amscd,amsfonts,amssymb,graphicx,color}

\firstpage{1} 
\pubvolume{10}
\issuenum{4}
\articlenumber{317}
\pubyear{2021}
\copyrightyear{2021}
\externaleditor{Academic Editor: Chris Goodrich} 
\datereceived{3 October 2021} 
\daterevised{17 November 2021} 
\dateaccepted{18 November 2021} 
\datepublished{24 November 2021} 
\hreflink{https://doi.org/10.3390/\linebreak axioms10040317}
\doinum{10.3390/axioms10040317}

\pdfoutput=1

\Title{Nabla Fractional Derivative and Fractional Integral on \mbox{Time Scales}}

\TitleCitation{Nabla Fractional Derivative and Fractional Integral on Time Scales}

\Author{Bikash Gogoi $^{1}$, Utpal Kumar Saha $^{1}$,  Bipan Hazarika $^{2},$ Delfim F. M. Torres $^{3}$\orcidA{} and  Hijaz Ahmad $^{4,}$*\orcidB{}}

\AuthorNames{Bikash Gogoi, Utpal Kumar Saha,  Bipan Hazarika, Delfim F. M. Torres and  Hijaz Ahmad}

\AuthorCitation{Gogoi, B.; Saha, U.K.; Hazarika, B.; Torres, D.F.M.; Ahmad, H.}

\address{%
$^{1}$ \quad Department of Basic and Applied Science, National Institute of Technology, 
Arunachal Pradesh, Jote 791113, India; bikash.phd20@nitap.ac.in (B.G.); uksahanitap@gmail.com (U.K.S.)\\
$^{2}$ \quad Department of Mathematics, Gauhati University, Guwahati 781014, Assam, India; bh\_rgu@yahoo.co.in\\
$^{3}$ \quad Center for Research and Development in Mathematics and Applications (CIDMA), 
Department of Mathematics, University of Aveiro, 3810-193 Aveiro, Portugal; delfim@ua.pt\\
$^{4}$ \quad Section of Mathematics, International Telematic University Uninettuno, 
Corso Vittorio Emanuele II, 39, \mbox{00186 Roma, Italy}}

\corres{Correspondence: hijaz555@gmail.com}

\abstract{In this paper, we introduce the nabla fractional derivative and fractional integral 
on time scales in the Riemann--Liouville sense. We also introduce the nabla fractional derivative 
in Gr\"unwald--Letnikov sense. Some of the basic properties and theorems related 
to nabla fractional calculus \mbox{are discussed.}}

\keyword{fractional calculus; Gr\"unwald--Letnikov derivative; 
Riemann--Liouville derivative; time scales calculus; nabla derivative; nabla integral\\
\\
\textbf{AMS:} Primary 26A33; Secondary 26E70}

\begin{document}
	

\section{Introduction}

Fractional calculus is a very significant branch in Mathematics whose applications are very useful for engineering 
students and researchers in both pure and applied field. The concept of ``fractional calculus'' was developed in 
the early of 17th century when L'Hopital asked Leibnitz the value of $\frac{1}{2}${th} order derivative. After that 
many mathematicians showed their interest on this topic. Initially the theory was developed mainly as a purely  
theoretical area. However, in the last decade it has been used in various fields such as mechanics, 
physics, chemistry, control theory, and many more, for instance one can see  \cite{LD,KE,MR,OS}.

The analysis of time scales calculus is a fairly new topic for researchers. Stefen Hilger and his Ph.D. supervisor, 
Bernd Aulbach, initiated the topic in the year 1988. After that, Hilger published two more paper on this topic 
\cite{HA,HD}. The theory was highly raised after the publication of two book on time scales by Martin Bohner 
and Allan Peterson \cite{MA,MD}. It combines the traditional areas of continuous and discrete analysis 
into one theory, which has various applications in discrete and continuous hybrid phenomena, 
quantum calculus and in various problems of economics \cite{FA}.

The inception of the idea of combining the time scales calculus and fractional calculus occurred 
in the Ph.D. dissertation of N.R.O. Bastos in 2012, where the delta (Hilger), and nabla derivative 
on time scales were discussed in fractional calculus using the tool of Laplace transform on some specific 
real and discrete time scales \cite{NDD,JL}. After the inception of the topic, a number of papers were 
published see \cite{ZL,NA,ES,GS,SM,DX,DY,ZY}. Recently, D. F. M. Torres introduced a generalized definition 
of Hilger derivative and integrals in a pure sense of Riemann--Liouville (RL) derivative \cite{NB,DC}. 
Many research works have been completed in a conformable delta and nabla fractional derivative 
and integrals \cite{BA,MM,Wangetal}.  


\section{Motivation of the Article} 

On the basis of above work, here we are motivated to study the nabla derivative and integral using 
a Gr\"unwald--Letnikov (GL) fractional derivative approach and then we arrive to the Riemann--Liouville 
sense. We introduce nabla fractional derivative and integral in unified approach of 
discrete and continuous time scales. Then, we generalize the definition of nabla fractional derivative 
and integral in arbitrary time scales and develop certain properties of nabla fractional 
derivative and fractional integral.

The paper is organized as follows. In Section~\ref{sec3}, we review briefly the essentials 
of time scales, as well as some basic definitions of nabla fractional derivative and integral 
which helps the readers to recognize easily our main findings. We assume that the readers are 
familiar with the basic view of time scales calculus and we refer the reader to go through 
\cite{MA,MD}. The paper also assess the Riemann--Liouville and Gr\"unwald--Letnikov fractional 
derivative and integral. Our main findings are given in Section~\ref{sec4} with some 
preliminaries definition and then we present fractional integral and fractional derivative 
in an arbitrary time scale $\mathbb{T}$. After that we prove certain important 
characteristics of fractional derivative and integral. We end with Section~\ref{sec5}
of conclusions.


\section{Preliminaries and Auxiliary Results}
\label{sec3}

\begin{Definition}[\cite{MA}]
\label{r1}
A time scale $\mathbb{T}$ is a closed subset of $\mathbb{R}$, 
with the subspace topology inherited from the stranded topology of $\mathbb{R}$.
The backward jump operator $\rho : \mathbb{T} \to \mathbb{T}$  is defined as  
$\rho(t) = \sup\left\{ s < t : s \in \mathbb{T} \right\}$ for $t \in \mathbb{T}$ 
and forward jump operator $\sigma(t)= \inf\left\{ s > t : s \in \mathbb{T}\right\}$. 
If $\rho(t) < t$ then $t$ is said to be a left scattered and if $\rho(t) = t$, 
then we say $t$ is a left dense point of $\mathbb{T}$, if $\sigma(t) > t$ and $\sigma(t) =t$, 
then we say $t$ is right scattered and right dense, respectively. Again, if $\mathbb{T}$ 
has a right scattered minimum $a$, then let $\mathbb{T}^{\kappa} = \mathbb{T} - \{a\}$, 
or else set $\mathbb{T}^{\kappa} = \mathbb{T}$. Here we consider the backward graininess 
$\nu : \mathbb{T}^{\kappa} \to [0, \infty]$, which is defined by $\nu(t) = t - \rho(t)$.
\end{Definition}   

\begin{Definition}[\cite{MD}]
\label{r2}
A function $h: \mathbb{T} \to \mathbb{R}$ is said to be a nabla differentiable at 
$t\in \mathbb{T}$, if for any $\varepsilon > 0$ there exists a neighborhood $V$ of $t$, 
such that 
$$
|h(\rho(t) - h(u) - h_{\nabla}(t)(\rho(t) - u)| \le \varepsilon|\rho(t) - u|
$$ 
for all $u\in V$. If $h_{\nabla}(t)$ exists for all $t \in \mathbb{T}^{\kappa}$ 
then it is called nabla derivative of $h$.
\end{Definition}

\begin{Theorem}[\cite{MD}]
\label{t1}
Let us consider a function $h : \mathbb{T} \to \mathbb{R}$ 
and let $t \in \mathbb{T}^{\kappa}$. Then we have
\begin{enumerate}
\item[(i)] If $h$ is continuous at a left-scattered $t$, then $h$ is nabla 
differentiable at $t$ with 
$$
h_{\nabla}(t) = \frac{h(t) - h(\rho(t))}{\nu(t)};
$$
\item[(ii)] If $t$ is left dense, then $h$ is nabla differentiable at $t$ 
if and only if the limit 
$$
\lim_{u \to t}\frac{h(t) - h(u)}{t - u}
$$ 
exists as a finite number. In this case 
$$
h_{\nabla}(t) = \lim_{u \to t}\frac{h(t) - h(u)}{t - u}.
$$
\end{enumerate}
\end{Theorem}

{\begin{Definition}[\cite{MD,JL}]
\label{r3} 
(Higher order nabla derivative): 
Assume a function $h : \mathbb{T} \to \mathbb{R}$, we first define the second order derivative 
$h_{\nabla \nabla}$ provided  $h_{\nabla}$ is differentiable on $\mathbb{T}^{K^{2}} 
= (\mathbb{T}^{\kappa})^{\kappa}$ with derivative $h_{\nabla \nabla} = (h_{\nabla})_{\nabla} 
: \mathbb{T}^{K^{2}} \to \mathbb{R}$. Similarly, proceeding up to $n^{th}$ order, here we 
obtain $h_{\nabla^{n}} : \mathbb{T}^{K^{n}} \to \mathbb{R}$, where  
$\mathbb{T}^{K^{n}}$ is a time scales which is obtained by removing $n$ right scattered left end point.
\end{Definition}}

\begin{Definition}[\cite{MR}]
\label{r4}
The Riemann--Liouville fractional differentiation of random order $\alpha$ 
is defined in the following manner:
$$
{\textsuperscript{RL}D^{-\alpha}_{x} 
= \frac{1}{\Gamma{(\alpha)}}
\int_{0}^{x}(x - t)^{\alpha - 1}h(t)dt}
$$  
for ${Re{(\alpha)} > 0}$. 

Riemann--Liouville derivative of order ${\alpha}\in{\mathbb{R}}$ is given by 
$$
{\textsuperscript{RL}D^{\alpha}_{x} 
= D^{n}\left\{\textsuperscript{RL}D^{\alpha - n}_{x}h(t)\right\} }
$$ 
for ${\alpha < n< \alpha + 1}$.
\end{Definition}

\begin{Definition}[\cite{KE}]
\label{r04} 
Let $\alpha > 0$. The Gr\"unwald--Letnikov derivative of fractional order $\alpha$ 
of a function $h$ is defined by 
$$
^{GL}D^{\alpha}h(t) = \frac{1}{h^{\alpha}}
\sum_{i = 0}^{\infty}{\alpha \choose i}h(t - ih), t \in \mathbb{R}.
$$ 
\end{Definition}

Referring Definition~\ref{r04} as the Gr\"unwald--Letnikov fractional derivative 
is quite common in literature (see \cite{AHM}). Moreover, once a starting point $t_{0}$ 
has been assigned, for practical reason then the following (see \cite{KN}) is often preferred,  
since it can be applied to function not defined (or simply not known) in $(-\infty, t_{0})$.

\begin{Definition}[\cite{BA}]
\label{r5}
A function $h : \mathbb{T} \to \mathbb{R}$ is called ld-continuous 
if $h$ is continuous at left dense point in $\mathbb{T}$ and right 
sided limit exists (finite) at right dense point in $\mathbb{T}$.
\end{Definition}

\begin{Definition}[\cite{DAH}]\label{r7}(Change of order of integration)
If we have any function $h(x, y, z)$ which is integrable with respect 
to $y$ and $z$, then the change of order of integration is given by the following formula:
$$
\int_{x_{0}}^{x}\int_{x_{0}}^{y}h(x, y, z)dzdy = \int_{x_{0}}^{x}\int_{z}^{x}h(x, y, z)dydz.
$$
\end{Definition}

\begin{Definition}[\cite{JL}]
\label{d1} 
Let $B\subset \mathbb{R}$. A subset $J\subset \mathbb{T}$ is called a time scales interval if 
$J = B\cap \mathbb{T}$. A function $ h : J \to \mathbb{T}$ is said to be left-dense 
absolutely continuous, if for all $\varepsilon > 0$ there exists $\delta > 0$, such that 
$\sum\limits_{k = 1}^{m}|h(y_{k}) - h(x_{k})|< \varepsilon$, whenever a disjoint finite 
collection of sub-time scales intervals $(x_{k}, y_{k})\cap \mathbb{T} \subset J$ 
for $1 \le k \le n$ satisfies $\sum\limits_{k = 1}^{m}|y_{k} - x_{k}| < \delta$. 
One denotes $h \in AC_{\nabla}$, if  $h_{\nabla}^{(n - 1)} \in AC$ 
then one denotes $h\in AC_{\nabla}^{(n)}$.
\end{Definition}


\section{Nabla Fractional Derivative and Nabla Fractional Integral}
\label{sec4}

\begin{Definition}
\label{r8}
For any time scale $\mathbb{T}$, a function $h : \mathbb{T} \to \mathbb{R}$ is said to be nabla
fractional differentiable of order $\mu$ at $ t \in \mathbb{T}^{\kappa}$, where 
$0 < \mu \le 1$, if for any $\varepsilon > 0$, there exists a neighborhood $V$ of $t$, such that 
\begin{equation}
\label{Eq:t1}
|h(\rho(t)) - h(u) - h^{(\mu)}_{\nabla}(t)(\rho(t) - u)^{\mu})| 
\le \varepsilon |(\rho(t)- u)^{\mu}|
\end{equation}
for all ${u \in V}$. If for all $t \in \mathbb{T}^{\kappa}$, $h$ holds the Equation (\ref{Eq:t1}), 
then we call ${h^{(\mu)}_{\nabla}(t)}$ the nabla fractional derivative of order ${\mu}$.
\end{Definition}

\begin{Theorem}
\label{r9} 
Nabla fractional derivative is not well defined in $\mathbb{T}$, 
but in $\mathbb{T}^{\kappa}$.
\end{Theorem}

\begin{proof} 
Let $h_{\nabla}^{(\mu)}(t)$ be defined at a point $t$ on a time scale $\mathbb{T}$, 
and assume that $t\notin \mathbb{T}^{\kappa}$. Then, $t \in \mathbb{T}\setminus 
\mathbb{T}^{\kappa}$. From Definition~\ref{r1}, $t$ must be unique which is equal 
to $a$, later, for any $\varepsilon > 0$ there exists a neighborhood 
$V =\left\{ t \right\}$ of $t$, we obtain for $u\in V$
$$
h(\rho(t)) = h(u) = h(\rho(a)) = h(a).
$$ 
Thus for $\zeta \in {\mathbb{R}}$ and $\mu \in (0, 1]$ we have
\begin{equation}
\label{e2} 
|h(\rho(t)) - h(u) - \zeta[(\rho(t) - u)^{\mu}]| 
= |h(a) - h(a) - \zeta(a - a)^{\mu}| \le  \varepsilon |(\rho(t) - u)^{\mu}|.
\end{equation}
Here, Equation (\ref{e2}) is true for each $\zeta \in \mathbb{R}$, 
which means for each $\zeta$ is the nabla derivative of $h$ of order $\mu$ if 
$t \notin \mathbb{T}^{\kappa}$, which cannot be true, so  $h_{\nabla}^{(\mu)}$ 
is well defined only on $\mathbb{T}^{\kappa}$.
\end{proof}

\begin{Theorem}
\label{r10}
For any time scale ${\mathbb{T}}$, let $h : \mathbb{T}^{\kappa} \to \mathbb{R}$. 
Then, for $\mu \in (0,1]$ we have the following:
\begin{enumerate}
\item[(i)] If $t$ is left dense and $h$ is nabla differentiable of order ${\mu}$ at $t$, 
then $h$ is continuous at $t$;
\item[(ii)] If $h$ is continuous at $t$ and $t$ is left scattered, 
then $h$ is nabla differentiable at $t$ of order $\mu$ with 
$$
{h^{(\mu)}_{\nabla}(t) = \frac{h(t) - h(\rho(t))}{\nu(t)^{\mu}}};
$$
\item[(iii)] If $t$ is left dense, then $h$ is differentiable at $t$ 
if and only if the limit 
$$
{\lim_{u \to t}\frac{h(t) - h(u)}{(t - u)^{\mu}}}
$$ 
exists as a finite number. In this case 
$$
{h^{(\mu)}_{\nabla}(t) = \lim_{u \to t}\frac{h(t) - h(u)}{(t - u)^{\mu}}};
$$
\item[(iv)] If $h$ is nabla differentiable of order $\mu $ at $t$, then 
$$
h(\rho(t)) = h(t) -  h^{(\mu)}_{\nabla}(t)(\nu(t))^{\mu}.
$$
\end{enumerate}
\end{Theorem}

\begin{proof}
(i) Given that $h$ is nabla fractional differentiable at $t$, 
then for $\varepsilon > 0$ there exists a neighborhood $V$ of $t$, 
such that 
$$
{|h(\rho(t)) - h(u) - h^{(\mu)}_{\nabla}(t)(\rho(t) - u)^{\mu}| 
\le \varepsilon |(\rho(t) - u)^{\mu}|}
$$
for $u \in V$. Therefore,
\begin{align*}
|h(t) - h(u)| \le 
&|h(\rho(t)) - h(u) - h^{(\mu)}_{\nabla}(t)(\rho(t) - u)^{\mu}| 
+ |h(\rho(t)) - h(t) - h^{(\mu)}_{\nabla}(t) (\rho(t) - t)^{\mu}|\\
+&|h^{(\mu)}_{\nabla}(t) ||(\rho(t) - u)^{\mu} + (\rho(t) - t)^{\mu}|
\end{align*}
for all $u \in V \cap (t - \varepsilon, t + \varepsilon)$ 
and since $t$ is a  left dense point, so
\begin{align*}
|h(t) - h(u)| 
&\le |h(\rho(t)) - h(u) - h^{(\mu)}_{\nabla}(t)(\rho(t) - u)^{\mu}| 
+ |h_{\nabla}^{(\mu)}(t)|(t - u)^{\mu}|\\
&\le \varepsilon |t - u|^{\mu} + |h_{\nabla}^{(\mu)}(t)|t - u|^{\mu}\\
& \le \varepsilon^{\mu}\big[ \varepsilon + |h_{\nabla}^{(\mu)}(t)|\big].
\end{align*}
It follows the continuity of $h$ at $t$.\\
(ii) Given that $h$ is continuous and $t$ is left scattered, 
by continuity 
$$
{\lim_{u \to t}\frac{h(\rho(t)) - h(u)}{(\rho(t) - u)^{\mu}} 
= \frac{h(\rho(t)) - h(t)}{(\rho(t) - t)^{\mu}}} 
= \frac{h(t) - h(\rho(t))}{(t - \rho(t))^{\mu}}.
$$
Hence, there exists a neighborhood $V$ of $t$, such that  
$$
{\left| \frac{h(\rho(t)) - h(u)}{(\rho(t) - u)^{\mu}} 
- \frac{h(t) - h(\rho(t))}{(t - \rho(t))^{\mu}} \right| 
\le \varepsilon}.
$$
For all $u \in V$, it follows that
$$
{\big| {h(\rho(t)) - h(u) - \frac{h(tt - h(\rho(t))}{(t - \rho(t))^{\mu}}(\rho(t) - u)^{\mu}} \big| 
\le \varepsilon|(\rho(t) - u)^{\mu}|}.
$$
From Definition~\ref{r8}, we obtain our result:
$$
{h^{(\mu)}_{\nabla}(t) = \frac{h(t) - h(\rho(t))}{\nu(t)^{\mu}}}.
$$
(iii) Given that $t$ is left dense, then we obtain $\rho(t) = t$, 
so there exists a neighborhood $V$ of $t$, such that 
$$
{|h(t) - h(u) - h^{(\mu)}_{\nabla}(t)(t - u)^{\mu}| \le \varepsilon |(t - u)^{\mu}|}
$$ 
for all $u \in V$. It follows that 
$$
\Big|\frac{h(t) - h(u)}{(t - u)^{\mu}} -h^{(\mu)}_{\nabla}(t)\Big| \le \varepsilon.
$$ 
So, we obtain 
$$
{h^{(\mu)}_{\nabla}(t) = \lim_{u \to t}\frac{h(t) - h(u)}{(t - u)^{\mu}}}.
$$
Now, assume that $\lim\limits_{u \to t}\frac{h(t) - h(u)}{(t - u)^{\mu}}$  
exists as a finite number, say $\mathcal{L}$, and $t$ is left dense. 
Then, for any $\varepsilon > 0$, 
$$
\left| \frac{h(t) - h(u)}{(t - u)^{\mu}} - \mathcal{L} \right| 
= \left| \frac{h(\rho(t) - h(u)}{(\rho(t) - u)^{\mu}} 
- \mathcal{L} \right| \le \varepsilon.
$$
Therefore, 
$$
|h(\rho(t)) - h(u) - \mathcal{L}(\rho(t) - u)^{\alpha}| \le \varepsilon|(\rho(t) - u)^{\mu}|
$$ 
from which we conclude that $h$ is fractional differentiable of order $\mu$ 
at $t$ and $h^{(\mu)}_{\nabla}(t) = \mathcal{L}$.\\
(iv) For all $t \in \mathbb{T}^{\kappa}$, there exist two possibilities of $t$. 

\noindent
{\bf Case 1:} If $t$ is left dense, then $\rho(t) = t$, 
and we have $\nu(t)^{\mu} = 0$, 
$$
h(\rho(t)) = h(t) -  h^{(\mu)}_{\nabla}(t)(\nu(t))^{\mu},
$$ 
so left hand side $\implies$ right hand side.

\noindent
{\bf Case 2:} If $t$ is left scattered, then $\rho(t) < t$, 
and by using Theorem~\ref{r10} (ii), we obtain  
$$
h(\rho(t)) = h(t) -  h^{(\mu)}_{\nabla}(t)(\nu(t))^{\mu}.
$$
The proof is complete.
\end{proof}

\begin{Proposition}
\label{p01} 
Let $h : \mathbb{T}^{\kappa} \to \mathbb{R}$. If $h(t) = k$ 
for all $t \in \mathbb{T}^{\kappa}$, and $k\in \mathbb{R}$, 
then ${h^{(\mu)}_{\nabla}(t) = 0}$ for all $\mu \in \mathbb{R}$.
\end{Proposition}

\begin{proof}
From the results (ii) and (iii) of Theorem~\ref{r10}, we have: 
if $t$ is left-scattered, then 
$$
{h^{(\mu)}_{\nabla}(t) = \frac{h(t) - h(\rho(t))}{(\nu(t))^{\mu}} 
= \frac{(k - k)}{(\nu(t))^{\mu}} = 0},
$$
if $t$ is left dense, then 
$$
h^{(\mu)}_{\nabla}(t) = \lim_{u \to t}\frac{h(t) -h(u)}{(t - u)^{\mu}} 
= \lim_{u \to t}\frac{(k - k)}{(t - u)^{\mu}} = 0,
$$ 
which completes the proof.
\end{proof}

\begin{Proposition}
\label{p2} 
For $h : \mathbb{T}^{\kappa} \to \mathbb{R}$, 
if $h(t) = t$, then for all $t\in \mathbb{T}^{\kappa}$ 
$$
h^{(\mu)}_{\nabla}(t) 
= \begin{cases}
\nu(t)^{1 - \mu}  
  & \mbox{~if~} \mu \neq 1,\\
1 & \mbox{~if~} \mu =1.
\end{cases}
$$
\end{Proposition}

\begin{proof}
Let $h_{\nabla}^{(\mu)}$ exist at $t \in \mathbb{T}^{\kappa}$. 
Then, from Theorem~\ref{r10} (iv), we have 
$$
{h(\rho(t)) = h(t) - h^{(\mu)}_{\nabla}(t)(\nu(t))^{\mu}},
$$ 
that is, ${\nu(t) = (\nu(t))^{\mu}h^{(\mu)}_{\nabla}(t)}$. 
If $\nu(t) \neq 0$, then we obtain our desired result which is 
$h^{(\mu)}_{\nabla}(t) = (\nu(t))^{1 - \mu}$ for $\mu \neq 1$.
If $\mu = 1$ then we obtain $h^{(\mu)}_{\nabla}(t) = (\nu(t))^{1 -1} = (\nu(t))^{0} =1$, 
hence the proof is complete.
\end{proof}

\begin{Theorem}
\label{t2}
Let $h,g : \mathbb{T} \to \mathbb{R}$ be two nabla differentiable 
functions of order $\mu \in (0, 1]$ at $t \in \mathbb{T}^{\kappa}$. 
Then, the following holds:
\begin{enumerate}
\item[(i)] The sum ${(\lambda_{1}h + \lambda_{2}g): \mathbb{T} \to \mathbb{R}}$ 
is nabla differentiable at $t$ of order $\mu$,  
$$
{(\lambda_{1}h + \lambda_{2}g)^{(\mu)}_{\nabla}(t) 
=\lambda_{1} h^{(\mu)}_{\nabla}(t) +\lambda_{2} g^{(\mu)}_{\nabla}(t)},
$$ 
where $\lambda_{1}$ and $\lambda_{2}$ are any two arbitrary constants;
\item[(ii)] The product $hg : \mathbb{T} \to \mathbb{R}$ is nabla 
differentiable of order $\mu$ at $t$,  
$$
(hg)^{(\mu)}_{\nabla} 
= h^{(\mu)}_{\nabla}(t)g(t) + h(\rho(t))g^{(\mu)}_{\nabla}(t)
 =  h(t)g^{(\mu)}_{\nabla}(t) + h^{(\mu)}_{\nabla}(t)g(\rho(t));
$$
\item[(iii)] If $h\left(t\right)h\left(\rho(t)\right) \neq 0$ 
then $\frac{1}{h}$ is nabla differentiable at $t$ of order $\mu$,
$$
\left(\frac{1}{h}\right)^{(\mu)}_{\nabla}(t) 
= - \frac{h^{(\mu)}_{\nabla}(t)}{h(t)h(\rho(t))};
$$
\item[(iv)] If $g(t)g(\rho(t)) \neq 0$, then $\frac{h}{g}$ 
is nabla differentiable at $t$ of order $\mu$,
$$
\left(\frac{h}{g}\right)_{\nabla}^{(\mu)}(t) 
= \frac{h^{(\mu)}_{\nabla}(t)g(t) - h(t)g^{(\mu)}_{\nabla}(t)}{g(t)g(\rho(t))}.
$$
\end{enumerate}
\end{Theorem}

\begin{proof}
(i) Let $\mu \in (0,1]$. Given that $h$ and $g$ are nabla differentiable 
at $t$ $\in \mathbb{T}^{\kappa}$ of order $\mu$, for any $\varepsilon > 0$ 
there exist neighborhoods $V_{1}$ and $V_{2}$ of $t$, thus for all $u \in V_{1}$
\begin{equation}
\label{Eq:e1}
|\lambda_{1}h(\rho(t)) - \lambda_{1}h(u) - \lambda_{1}h^{(\mu)}_{\nabla}(t)((\rho(t) - u)^{\mu}| 
\le \frac{\varepsilon}{2}|(\rho(t) - u)^{\mu}|,
\end{equation} 
also for all $u\in V_{2}$
\begin{equation}
\label{Eq:e2}
|\lambda_{2}g(\rho(t)) - \lambda_{2}g(u) 
- \lambda_{2}g^{(\mu)}_{\nabla}(t)((\rho(t) - u)^{\mu}| 
\le \frac{\varepsilon}{2}|(\rho(t) - u)^{\mu}|.
\end{equation} 
Let $u \in V = V_{1} \cap V_{2}$. Then, we obtain
\begin{align}
\label{e3}
\begin{split}
& \quad |(\lambda_{1}h + \lambda_{2}g)(\rho(t)) - (\lambda_{1}h + \lambda_{2}g)(u) 
- [\lambda_{1}h^{(\mu)}_{\nabla}(t) + \lambda_{2}g^{(\mu)}_{\nabla}(t)]((\rho(t) - u)^{\mu}|\\ 
&\le|[\lambda_{1}h(\rho(t)) - \lambda_{1}h(u) 
- \lambda_{1}h^{(\mu)}_{\nabla}(t)((\rho(t) - u)^{\mu}]|\\
&+ |[\lambda_{2}g(\rho(t)) - \lambda_{2}g(u) 
- \lambda_{2}g^{(\mu)}_{\nabla}(t)((\rho(t) - u)^{\mu}]|.
\end{split}
\end{align}
By using the Equations (\ref{Eq:e1}) and (\ref{Eq:e2}), we obtain that
\begin{align*}
\begin{split}
&\quad |(\lambda_{1}h + \lambda_{2}g)(\rho(t)) - (\lambda_{1}h + \lambda_{2}g)(u) 
- [\lambda_{1}h^{(\mu)}_{\nabla}(t) + \lambda_{2}g^{(\mu)}_{\nabla}(t)]((\rho(t) - u)^{\mu}|\\
&< \frac{\varepsilon}{2}|(\rho(t) - u)^{\mu}| + \frac{\varepsilon}{2}|(\rho(t) - u)^{\mu}|\\ 
&=\varepsilon |(\rho(t) - u)^{\mu}|.
\end{split}
\end{align*}
From Theorem~\ref{r10}, it holds that 
$(\lambda_{1}h + \lambda_{2}g)^{(\mu)}_{\nabla}$ is a nabla differentiable 
at $t \in \mathbb{T}^{\kappa}$ of order $\mu$.\\
(ii) If $t$ is left dense, i.e., $\rho(t) = t$ for $t \in \mathbb{T}^{\kappa}$, then 
\begin{align*}
(hg)^{(\mu)}_{\nabla}(t)
&=\lim_{u \to t}\frac{(hg)(t) - (hg)(t)}{(t - u)^{(\mu)}}\\
&= \lim_{u \to t}\frac{h(t) - h(u)}{(t - u)^{\mu}}g(t) 
+ \lim_{u \to t}\frac{g(t) - g(u)}{(t - u)^{\mu}}h(u)\\
&= h^{(\mu)}_{\nabla}(t)g(t) + g^{(\mu)}_{\nabla}(t)h(t)\\
&=  h^{(\mu)}_{\nabla}(t)g(t) + g^{(\mu)}_{\nabla}(t)h(\rho(t)).
\end{align*}
Additionally, if $\rho(t) < t$, then
\begin{align*}
(hg)^{(\mu)}_{\nabla}(t) 
&= \frac{(hg)(t) - (hg)(\rho(t))}{(\nu(t))^{\mu}}\\
&= \frac{h(t) - h(\rho(t))}{(\nu(t))^{\mu}}g(t) 
+ \frac{g(t) - g(\rho(t))}{(\nu(t))^{\mu}}f(\rho(t))\\
&= h^{(\mu)}_{\nabla}(t)g(t) + h(\rho(t))g^{(\mu)}_{\nabla}(t).
\end{align*}
Other part of the proof is very similar to this.\\
(iii) Using the above result and Proposition~\ref{p01}, we obtain 
\begin{align*} 
\left(h \cdot \frac{1}{h}\right)^{(\mu)}_{\nabla}(t)= (1)^{(\mu)}_{\nabla}(t)&=0
\end{align*}
and hence, by (ii),
\begin{align*}
\left(\frac{1}{h}\right)^{(\mu)}_{\nabla}(t)h{(\rho(t))} 
+ h^{(\mu)}_{\nabla}(t)\frac{1}{h(t)} &= 0.
\end{align*}
Since $h(t)h(\rho(t)) \ne 0$, so we obtain
\begin{align*}
\left(\frac{1}{h}\right)^{(\mu)}_{\nabla}(t)
&= -\frac{h^{(\mu)}_{\nabla}(t)}{h(t)h(\rho(t))}.
\end{align*}
(iv) Using the result of Theorem~\ref{t2} (ii) and (iii) we obtain the following:
\begin{align*}
\left(\frac{h}{g}\right)^{(\mu)}_{\nabla}(t) 
&= \left(h \cdot \frac{1}{g}\right)^{(\mu)}_{\nabla}(t)\\
& = h(t)\left(\frac{1}{g}\right)^{(\mu)}_{\nabla}(t) + h^{(\mu)}_{\nabla}(t)\frac{1}{g(\rho(t))}\\
& = -h(t)\frac{g^{(\mu)}_{\nabla}(t)}{g(t)g(\rho(t))} + h^{(\mu)}_{\nabla}(t)\frac{1}{g(\rho(t))}\\
&= \frac{ h^{(\mu)}_{\nabla}(t)g(t) - h(v)g^{(\mu)}_{\nabla}(t)}{g(t)g(\rho(v))}. 
\end{align*}
This completes the proof.
\end{proof}

\begin{Theorem}
\label{t3}
Let $k$ be a constant, $n \in \mathbb{N}$. Then, for $0 < \mu \le 1$, 
we obtain the following:
\begin{enumerate} 
\item[(i)] If $h(t) = (t - k)^{n}$, then
$$
h^{\mu}_{\nabla}(t) = (\nu(t))^{1 - \mu}\sum^{n - 1}_{j = 0}(\rho(t) - k)^{j}(t - k)^{n - 1 -j};
$$
\item[(ii)] If $g(t) = \frac{1}{(t - k)^{n}}$, then
$$
g^{(\mu)}_{\nabla}(t) = -(\nu(t))^{1 - \mu}\sum^{n - 1}_{j = 0}
\frac{1}{(\rho(t) - k)^{n - j}(t - k)^{j + 1}}
$$ 
provided $(t - k)(\rho(t) - k) \neq 0$.
\end{enumerate}
\end{Theorem}

\begin{proof}
(i) Here we prove this result by using the method of induction. 
If $n = 1$, then $h(t) = t - k$ hence 
$h^{(\mu)}_{\nabla}(t)  = (\nu(t))^{1 - \mu}$ 
is true from Propositions~\ref{p1} and \ref{p2}. We assume that 
$$
h^{\mu}_{\nabla}(t) = (\nu(v))^{1 - \mu}\sum^{n - 1}_{j = 0}(\rho(t) - k)^{j}(t - k)^{n - 1 -j}
$$ 
holds for $h(t) = (t - k)^{n}$. We shall prove the result is true for 
$$
H(t) = (t - k)^{n + 1} = (t - k)h(t). 
$$
By using Theorem~\ref{t2} (ii), we obtain 
\begin{align*}
h^{(\mu)}_{\nabla}(t) 
&= (t - k)^{(\mu)}_{\nabla}h(\rho(t)) + h^{(\mu)}_{\nabla}(t)(t - k) 
= (\nu(t))^{1 - \mu}h(\rho(t)) + h^{(\mu)}_{\nabla}(t)(t - k)\\
&=(\nu(t))^{1 - \mu}(\rho(t) - k)^{n} + (\nu(t))^{1 - \mu}(t - k)
\sum^{n - 1}_{j = 0}(\rho(t) - k)^{j}(t - k)^{n - 1 - j}\\
&=(\nu(t))^{1 - \mu}\big[(\rho(t) - k)^{n} + \sum^{n - 1}_{j = 0}(
\rho(t) - k)^{j}(t - k)^{n - j}\big]\\
&=(\nu(t))^{1 - \mu} \sum^{n }_{j = 0}(\rho(t) - k)^{j}(t - k)^{n - j}.
\end{align*} 
(ii) Let $g(t) = \frac{1}{(t - k)^{n}} = \frac{1}{h(t)}$. 
Using Theorem~\ref{t2} (iii), we obtain
\begin{align*} 
g^{(\mu)}_{\nabla}(t)
&= \left(\frac{1}{h(t)}\right)^{(\mu)}_{\nabla}(t) \\
&= -\frac{h^{(\mu)}_{\nabla}(t)}{h(t)h(\rho(t))} \\
&= -(\nu(t))^{1 - \mu}\frac{\sum\limits_{j = 0}^{n - 1}(\rho(t) 
- k)^{j}(t - k)^{n - 1 - j}}{(t - k)^{n}(\rho(t) - k)^{n}}\\
&= -(\nu(t))^{1 - \mu}\sum^{n - 1}_{j = 0}\frac{1}{(t - k)^{j + 1}(\rho(t) - k)^{n - j}},
\end{align*}
provided $(t - k)(\rho(t) - k) \neq 0$.
\end{proof}

\begin{Example}
\label{ex1}
Let $\mu \in (0,1]$.
\begin{enumerate}
\item[(i)] If $g(t) = t^{2}$, then from Theorem~\ref{t3}, we obtain
\begin{align*}
g^{(\mu)}_{\nabla}(t)
&= (\nu(t))^{1 - \mu}\left[\sum^{1}_{j = 0}(\rho(t))^{j}(t)^{1 - j}\right]\\
&= (\nu(t))^{1 - \mu}\left[(\rho(t))^{0}(t)^{1} + (\rho(t))^{1}(t)^{0}\right]\\
&= (\nu(t))^{1 - \mu}(t + \rho(t)).
\end{align*}
By using  Theorem~\ref{t3}, we obtain the following results:
\item[(ii)] If $g(t) = t^{3}$, then  
$g^{(\mu)}_{\nabla}(t) = (\nu(t))^{1 - \mu}[ t^{2} + t\rho(t)  + (\rho(t))^{2}]$.
\item[(iii)] If $g(t) = \frac{1}{t}$, 
then $g^{(\mu)}_{\nabla}(t) = -\frac{(\nu(t))^{1 - \mu}}{t(\rho(t))}$.
\item[(iv)] If $g(t) = \frac{1}{t^{2}}$, 
then $g_{\nabla}^{(\mu)}(t) = -\frac{\nu(t)^{1 - \mu}}{t(\rho(t))^{2} + t^{2}\rho(t)}$.
\end{enumerate}
\end{Example}

\begin{Corollary}
\label{cor1}
Nabla fractional derivative in some specific time scales $\mathbb{T}$.
\begin{enumerate}
\item[(i)] If we consider the real time scale $\mathbb{T} = \mathbb{R}$, 
then all the elements of $\mathbb{T}$ are dense. So, by using  
Theorem~\ref{r10} (iii), we have that
$$ 
h^{(\mu)}_{\nabla}(t) = \lim_{u \to t}\frac{h(t) - h(u)}{(t - u)^{\mu}}
$$ 
exists, if $\mu =1$, then we have  $h_{\nabla}^{(\mu)} = h^{\prime}(t)$,
which is similar to the ordinary derivative. 
\item[(ii)] If $\mathbb{T} = \mathbb{Z}$, for $t \in \mathbb{T}$ one has 
$\rho(t) = t - 1$ and then $\nu(t) = t - (t - 1) = 1$. Now, 
by using Theorem~\ref{r10} (ii), we obtain
\begin{align*}
h^{(\mu)}_{\nabla}(t) 
&= \frac{h(t) - h(\rho(t))}{(t - \rho(t))^{\mu}}\\ 
&= \frac{h(t) - h(t - 1)}{(t  - (t - 1))^{\mu}}.
\end{align*} 
If $\alpha = 1$ we have that $h^{(\mu)}_{\nabla}(t) = h(t) - h(t - 1) = \nabla{h(t)}$,  
which is similar as the usual backward operator;
\item[(iii)] Let $\mathbb{T} = h\mathbb{Z}$, where $h > 0$. Then we obtain
$$ 
\rho(t) = \sup\left\{u \in \mathbb{T}: u < t \right\} 
= \sup\left\{t - nh: \mbox{~for~} n  \in \mathbb{N}\right\} = t - h
$$
and then the function $\nu(t) = t - \rho(t) = t - (t - h) =  h$, which is constant.\\
For $g : h\mathbb{Z} \to \mathbb{R}$, we have from Theorem~\ref{t1} that
\begin{align*}
g_{\nabla}(t) = g^{(1)}_{\nabla}(t) = \frac{g(t - h(\rho(t))}{\rho(t) - t} 
= \frac{g(t - h) - g(t)}{(t - h - t)} = \frac{g(t) - g(t - h)}{h}.
\end{align*}
From Definition~\ref{r3}, the second order nabla derivative is  
\begin{align*}
g_{\nabla \nabla}(t) = g^{(2)}_{\nabla}(t) 
= \frac{g_{\nabla}(t) - g_{\nabla}(t - h)}{h} 
= \frac{g(t) - 2g(t - h) + g(t - 2h)}{h^{2}}.
\end{align*}
In general, the $m^{th}$ derivative for 
$t \in h\mathbb{Z}$ and $m \in \mathbb{N}$,
$$
g^{(m)}_{\nabla}(t) = \frac{1}{h^{m}}\sum^{m}_{r = 0}(-1)^{r}{m \choose j}g(t - rh),
$$ 
where the binomial coefficient ${m \choose r}$ is defined as follows: 
\begin{equation}
\label{Eq:4}
\begin{split}
{m \choose r} &= \frac{m(m - 1)(m - 2)\cdot\cdot\cdot(m - r + 1)}{r!}\\
&= \begin{cases}
\frac{m!}{r!(m - r)!}, & r = 0,1,2,\cdot\cdot\cdot\\
0, & r > m.
\end{cases}
\end{split}
\end{equation}
Since the binomial coefficient vanish when $r > m$, so no contribution 
in the summation is given from the presence of terms with $r > m$, 
the upper limit of the formula can be raised to any value greater 
than $m$ and hence, the finite summation in this formula can be replaced 
with the infinite series, i.e.,
$$
g^{(m)}_{\nabla}(t) = \frac{1}{h^{m}}\sum^{\infty}_{r = 0}(- 1)^{r}{m \choose j}g(t - rh).
$$
Letting $h$ tend to zero, then all points of the time scale become dense, 
and the time scale becomes the continuous time scale. If the value of $m$ 
is replaced by an arbitrary real number $\mu \in \mathbb{R}$, $\mu > 0$, 
and changing the factorial function with a Euler gamma function using 
the recurrence relation $(n - 1)! = \Gamma(n)$, then without losing 
the generality, if we replace $m$ by any arbitrary real number 
$\mu \in \mathbb{R}$, then the nabla fractional derivative, 
from Definition~\ref{r3} and Theorem \ref{r10}, is
\begin{equation}\label{Eq:5}
^\mathbb{T}D^{\mu}_{t} = g^{(\mu)}_{\nabla}(t) 
= \lim_{h \to 0}\frac{1}{h^{\mu}}\sum^{\infty}_{r = 0}(- 1)^{r}{\mu \choose r}g(t - rh).
\end{equation}
Moreover, once a starting point $a$ assign as $nh = t - a$ for $t > a$, such that
\begin{equation}
\label{Eq:6}
_{a}^{\mathbb{T}}D^{\mu}_{t}g(t) = \lim_{h \to 0}\frac{1}{h^{\mu}}
\sum_{r = 0}^{n}(-1)^{r}{\mu \choose r}g(t - rh). 
\end{equation}
Since for any continuous function $g(t)$ 
Gr\"unwald--Letnikov derivative and Riemann--Liouville 
derivative coincide with positive non integer order derivative, so we have 
$$
_{a}^{\mathbb{T}}D^{\mu}_{t}g(t) = ^{RL}_{a}D^{\mu}_{t}g(t),
$$ 
where $^{RL}_{a}D^{\mu}_{t}g(t)$ denote the Riemann--Liouville 
fractional derivative defined on time scales, which is most useful 
in the study of fractional calculus.\\
If $\mu < 0$, then we have
\begin{equation}
\label{Eq:7}
\begin{bmatrix} 
\mu \\ 
r 
\end{bmatrix} 
= \frac{(\mu)(\mu + 1)(\mu + 2)(\mu + 3)\cdots (\mu + r - 1)}{r!},
\end{equation}
i.e., when $\mu = -\mu$, then from Equations (\ref{Eq:4}) and (\ref{Eq:7}) 
we obtain
$$
{-\mu \choose r} = \frac{(-\mu - 1)(-\mu - 2)(-\mu - 3)
\cdots (-\mu - r + 1)}{r!}
= (-1)^{r}
\begin{bmatrix} 
\mu \\ 
r 
\end{bmatrix}
$$
or we can write $ (-1)^{r}{\mu \choose r} 
= \begin{bmatrix} \mu \\ r \end{bmatrix}$. 
For any ld-continuous function and for $\mu = - \mu$, 
then from Equation (\ref{Eq:6}) we obtain 
\begin{align}
\label{Eq:8}
\begin{split}
_{a}^{\mathbb{T}}D^{-\mu}_{t}g(t)
&=  \lim_{h \to 0}h^{\mu}\sum_{r = 0}^{n}
\begin{bmatrix}
\mu \\  r
\end{bmatrix} 
g(t - rh) \\
&=\lim_{h \to 0}\sum_{r = 0}^{n}\frac{1}{r^{\mu - 1}} 
\begin{bmatrix} 
\mu \\ 
r 
\end{bmatrix}
h(rh)^{\mu - 1}g(t - rh)\\
& = \frac{1}{\Gamma(\mu)}\lim_{h \to 0}\sum_{h \to 0}^{n}\frac{\Gamma(\mu)}{r^{\mu - 1}} 
\begin{bmatrix} 
\mu \\ 
r 
\end{bmatrix}
h(rh)^{\mu - 1}g(t - rh)\\
&=\frac{1}{\Gamma(\mu)}\lim_{n \to \infty}
\sum_{r = 0}^{n}\frac{\Gamma(\mu)}{r^{\mu - 1}}
\begin{bmatrix} 
\mu \\ 
r
\end{bmatrix} 
\Big( \frac{t - a}{n}\Big)\Big( r\frac{t - a}{n} \Big)^{\mu - 1}
g\Big(t - r\frac{t - a}{n}\Big).
\end{split}
\end{align}
Let us take 
$$
\beta_{r} = \frac{\Gamma(\mu)}{r^{\mu - 1}}
\begin{bmatrix} 
\mu \\ 
r 
\end{bmatrix}, 
\eta_{n,r} 
= \Big( \frac{t - a}{n}\Big)\Big( r\frac{t - a}{n} \Big)^{\mu - 1}
g\Big( t - r\frac{t - a}{n}\Big).
$$ 
Now,
\begin{equation}
\label{Eq:9}
\lim_{r \to \infty}\beta_{r} 
= \lim_{r \to \infty}\frac{\Gamma(\mu)}{r^{\mu - 1}}
\begin{bmatrix} 
\mu \\ 
r 
\end{bmatrix} = 1.
\end{equation}
\begin{align}
\label{Eq:10}
\begin{split} 
\lim_{n \to \infty}\sum_{r = 0}^{n}\eta_{n,r} 
&=\lim_{n \to \infty}\sum_{r = 0}^{n} \Big( \frac{t - a}{n}\Big)\Big( 
r\frac{t - a}{n} \Big)^{\mu - 1}g\Big( t - r\frac{t - a}{n}\Big)\\
&= \lim_{n \to \infty}\sum_{r = 0}^{n}h(rh)^{\mu - 1}g(t - rh)\\
&= \int_{a}^{t}\big( t - s \big)^{\mu - 1}g(s)\nabla s.
\end{split}
\end{align}
Here, we obtain a condition  (see \cite{IF}) that, if 
\begin{equation}
\label{Eq:11}
\lim_{r \to \infty}\beta_{r} = A 
\text{ and } \lim_{n \to \infty}
\sum_{r = 0}^{n}\eta_{n,r} = B, 
\text{ then } 
\lim_{n \to \infty}\sum_{r = 0}^{n}\eta_{n,r}\beta_{r} = AB.
\end{equation}
By using the Equations (\ref{Eq:8})--(\ref{Eq:11}), we obtain
\begin{equation}
\label{Eq:12}
_{a}^{\mathbb{T}}D^{-\mu}_{t}g(t) 
= \frac{1}{\Gamma(\mu)}
\int_{a}^{t}(t - s)^{\mu - 1}g(s)\nabla s, 
\end{equation}
which represents the nabla integral of any arbitrary 
order $\mu$ in a Riemann--Liouville sense.
\end{enumerate}
\end{Corollary} 

\begin{Remark}
\label{rem1}
The definition of nabla fractional integral defined in 
Equation (\ref{Eq:12}) is not the natural one 
for arbitrary time scales $\mathbb{T}$. For showing this 
we take an example. If $g(t) = t^{2}$, then from 
Example~\ref{ex1}, $g^{\mu}_{\nabla}(t) = t + \rho(t)$, 
for $\mu = 1$. If the time scale is the continuous time scale 
$\mathbb{T} = \mathbb{R}$, then $\rho(t) = t$ and, hence, 
from Corollary~\ref{cor1}, we find that 
$g_{\nabla}^{(\mu)}(t) = g^{\prime}(v) = 2t$. 
But if we take the discrete time scale, we obtain $\rho(t) = t - 1$, 
and the nabla derivative on $\mathbb{T} = \mathbb{Z}$ means 
the backward difference of $t^{2}$, i.e., $\nabla(t^{2}) = g(t) - g(t - 1) = 2t - 1$. 
Again, since every ld-continuous function is nabla integrable, 
so in this case we can claim that $\int_{0}^{t}g_{\nabla}^{(\mu)} = t^{2}$. 
In $\mathbb{T} = \mathbb{R}$ also means that $\int_{a}^{b}f(s)\nabla s = \int_{a}^{b}f(s)ds$, 
so $\int_{0}^{t}(s + \rho(s))\nabla s = \int_{0}^{t}(2s)ds =  t^{2}$, 
but in $\mathbb{T} = \mathbb{Z}$ we obtain $\int_{a}^{b}f(s)\nabla s 
= \sum\limits_{a}^{b}f(s)$ for $a < b$, so 
$\int_{0}^{t}(s + \rho(s))\nabla s = \sum\limits_{0}^{t}(2s - 1) = t^{2}$ for $t > 0$. 
For the appearance of $\rho(s)$, we claim that the correct definition 
of nabla fractional integral on an arbitrary time scale $\mathbb{T}$ must be
\begin{equation}
\label{Eq:13} 
_{a}^{\mathbb{T}}\mathbb{I}^{\mu}_{t}g(t) 
= \frac{1}{\Gamma(\mu)}\int_{a}^{t}(t - \rho(s))^{\mu - 1}g(s)\nabla s,
 \end{equation}
which is the generalization of the nabla fractional integral defined 
on the Equation (\ref{Eq:12}), in a Riemann--Liouville sense.
\end{Remark}

\begin{Proposition}
\label{p02} 
The nabla fractional integral for any function $g$ defined on $[a,b]$ satisfies
$$
_{a}^{\mathbb{T}}\mathbb{I}_{t}^{\mu}\circ~ _{a}^{\mathbb{T}}\mathbb{I}_{t}^{\beta} 
=~ _{a}^{\mathbb{T}}\mathbb{I}_{t}^{\mu + \beta}
$$ 
for $\mu > 0, \beta > 0$.
\end{Proposition}

\begin{proof} 
By using the generalized definition of nabla derivative 
of fractional order from the Equation (\ref{Eq:13}), we have
\begin{align*}
\Big({_{a}^{\mathbb{T}}\mathbb{I}_{t}^{\mu}
\circ~ _{a}^{\mathbb{T}}\mathbb{I}_{t}^{\beta}}\Big)(g(t)) 
&= ~ _{a}^{\mathbb{T}}\mathbb{I}_{t}^{\mu}\Big( { _{a}^{\mathbb{T}}\mathbb{I}_{t}^{\beta}}(g(t))\Big)\\
&= \frac{1}{\Gamma(\mu)}\int_{a}^{t}\big({t - \rho(s)}\Big)^{\mu 
- 1}\Big({_{a}^{\mathbb{T}}\mathbb{I}_{t}^{\beta}(g(s))}\Big) \nabla s\\
&= \frac{1}{\Gamma(\mu)}\int_{a}^{t}\big( t- \rho(s))^{\mu - 1}\frac{1}{\Gamma(\beta)}
\int_{a}^{s}( s - \rho(v))^{\beta - 1}g(v)\nabla v)\Big) \nabla s\\
&= \frac{1}{\Gamma(\mu)\Gamma(\beta)}\int_{a}^{t}\int_{a}^{s}\Big( 
{(t - \rho(s))^{\mu - 1}( s - \rho(v))^{\beta - 1}g(v)\nabla v}\Big)\nabla s.
\end{align*}
By using Definition~\ref{r7}, we obtain
\begin{align*}
\Big({_{a}^{\mathbb{T}}\mathbb{I}_{t}^{\mu}
\circ~ _{a}^{\mathbb{T}}\mathbb{I}_{t}^{\beta}}\Big)(g(t)) 
= ~ \frac{1}{\Gamma(\mu)\Gamma(\beta)}\int_{a}^{t} \Big[  
\int_{v}^{t}{(t - \rho(s))^{\mu - 1}( s - \rho(v))^{\beta - 1}
\nabla s}\Big]g(v)\nabla u.
\end{align*}
Let $s = \rho(v) + x(t - \rho(v))$ for $x \in \mathbb{R}$. We have 
\begin{align*}
&\quad \Big({_{a}^{\mathbb{T}}\mathbb{I}_{t}^{\mu}
\circ~ _{a}^{\mathbb{T}}\mathbb{I}_{t}^{\beta}}\Big)(g(t))\\
&=~  \frac{1}{\Gamma(\mu)\Gamma(\beta)}\int_{a}^{t} \Big[ 
{\int_{0}^{1}(1 - x)^{\mu - 1}(t - \rho(v))^{\mu - 1}x^{\beta - 1}
(t - \rho(v))^{\beta - 1}(t - \rho(v))dx} \Big] g(v)\nabla v\\
&=  \frac{1}{\Gamma(\mu)\Gamma(\beta)}\int_{0}^{1}(1 - x)^{\mu - 1}
x^{\beta - 1}dx\int_{a}^{t}(t - \rho(v))^{\mu + \beta - 1}g(v)\nabla v\\
&= \frac{B(\mu,\beta)}{\Gamma(\mu)\Gamma(\beta)}\int_{a}^{t}(t - \rho(v))^{\mu + \beta - 1}
g(v)\nabla v\\
&= \frac{1}{\Gamma(\mu + \beta)}\int_{a}^{t}(t - \rho(v))^{\mu + \beta - 1}g(v)\nabla v\\
&=~ _{a}^{\mathbb{T}}\mathbb{I}_{t}^{\mu + \beta}g(t).
\end{align*}
This completes the proof.
\end{proof}

Next definition uses integration as an anti-derivative process. 

\begin{Definition}
\label{r21}
(Riemann--Liouville fractional derivative on time scales)  
For $t\in \mathbb{T}$ and $g : \mathbb{T} \to \mathbb{R}$, 
the (left) Riemann--Liouville fractional derivative  
of order $\mu \in (0,1]$ is defined by 
\begin{align*}
_{a}^{RL}D_{t}^{\mu}g(t) 
= \frac{1}{\Gamma(1 - \mu)}\int_{a}^{t} 
{\Big( {t - \rho(s))^{-\mu}g(s)\nabla s } \Big)^{\nabla}}.
\end{align*}
\end{Definition}

\begin{Remark}
\label{r22} 
If $\mathbb{T} = \mathbb{R}$, then Definition~\ref{r21} 
gives the classical (left) Riemann--Liouville derivative  
of fractional order $\mu$.  Here, we are only studying 
the derivative in terms of left operators, the analogous 
right operators are easily acquired by changing the limit 
of integration. 
\end{Remark}

A different extension to time scales is obtained
by using the nabla fractional derivative in terms of Caputo 
sense, that will be more effective for integer order initial 
conditions and are more easy to obtain in real world problems
\cite{MR,JL}.

\begin{Definition}
\label{r23}
(Nabla derivative on time scales in a Caputo sense)
For {$t, t_{0} \in \mathbb{T}$, let us assume a finite time scale 
interval $[t_{0}, t]\cap \mathbb{T}^{{\kappa}^{n}} 
= [t_{0}, t]_{\mathbb{T}^{{\kappa}^{n}}}$. Then, 
for any  $g \in AC_{\nabla}^{(n)}[t_{0}, t]_{\mathbb{T}^{{\kappa}^{n}}}$ 
of absolutely continuous function as in Definition~\ref{d1}}, 
we define the Caputo nabla fractional derivative of order $0 < \mu \le 1$ as
$$
_{a}^{C}D_{t}^{\mu}g(t) 
= \frac{1}{\Gamma(n - \mu)}
\int_{a}^{t}(t - \rho(s))^{n -\mu - 1}
g_{\nabla}^{(n)}(s)\nabla s,
$$ 
where $n = [\mu]+ 1$.
\end{Definition} 

\begin{Definition}
\label{r24}
If $\mu < 0$, then the nabla derivative of order $\mu$ 
in terms of Riemann--Liouville, is the fractional integral 
of order $-\mu$, that is,  $_{a}^{RL}D_{t}^{\mu} 
=  _{a}^{\mathbb{T}}\mathbb{I}^{-\mu}_{t}$.
\end{Definition}

\begin{Definition}
\label{r25} 
If $\mu < 0$, then the nabla fractional integral of order $\mu$ 
is the nabla fractional derivative of order $-\mu$, i.e.,  
$_{a}^{\mathbb{T}}\mathbb{I}_{t}^{\mu} = ~ _{a}^{RL}D^{-\mu}_{t}$.
\end{Definition}

\begin{Proposition}
\label{p1}
Let $g : \mathbb{T}^{\kappa} \to \mathbb{R}$ be a nabla fractional 
differentiable function. Then, for any $0 < \mu \le 1$, 
$$ 
_{a}^{RL}D_{t}^{\mu}g(t) = 
\nabla \circ~{ _{a}^{\mathbb{T}}\mathbb{I}_{t}^{1 - \mu}}g(t).
$$ 
\end{Proposition}

\begin{proof} 
Let $g : \mathbb{T}^{\kappa} \to \mathbb{R}$ be a nabla fractional 
differentiable function. Then, from Definition~\ref{r21} and 
Equation (\ref{Eq:13}), we have 
$$
^{RL}_{a}D^{\alpha}_{t}g(t) 
=\frac{1}{\Gamma(1 - \alpha)}
\int^{t}_{a}\left((t - \rho(t))^{-\alpha}g(s)\nabla s\right)^{\nabla} 
= \big({^{\mathbb{T}}_{a}\mathbb{I}^{1 -\alpha}_{t}g(t)}\big) 
= \big({\nabla \circ~ _{a}^{\mathbb{T}}\mathbb{I}^{1 - \alpha}_{t}}\big)g(t).
$$
The proof is complete.
\end{proof}

\begin{Proposition}
\label{p2:b} 
For any integrable function $g$ defined on a time scales 
interval $[a,b]$ one has
\begin{align*}
_{a}^{RL}D_{t}^{\mu} 
\circ~ _{a}^{RL}\mathbb{I}_{t}^{\mu}g(t) =  g(t).
\end{align*}
\end{Proposition}

\begin{proof} 
From Propositions~\ref{p01} and \ref{p1}, we obtain
\begin{align*}
_{a}^{RL}D_{t}^{\mu} 
\circ~ _{a}^{\mathbb{T}}\mathbb{I}_{t}^{\mu}g(t)
& =  \Big[ { _{a}^{\mathbb{T}}\mathbb{I}_{t}^{1 - \mu}\Big({
_{a}^{\mathbb{T}}\mathbb{I}_{t}^{\mu}(g(t))\Big)}}\Big]^{\nabla}\\
& = \Big[{_{a}^{\mathbb{T}}\mathbb{I}_{t}g(t)}\Big]^{\nabla} = g(t).
\end{align*}
This concludes the proof.
\end{proof}

\begin{Corollary}
\label{cor2}
For $0 < \mu \le 1$, we have  
$ _{a}^{RL}D_{t}^{\mu}\circ~ _{a}^{RL}D_{t}^{-\mu} = I$ 
and $_{a}^{\mathbb{T}}\mathbb{I}_{t}^{-\mu} 
\circ{ _{a}^{\mathbb{T}}\mathbb{I}_{t}^{\mu}} = I$, 
where $I$ denotes the identity operator.
\end{Corollary}

\begin{proof}
From Definition~\ref{r24} and Proposition~\ref{p2:b}, we find that  
$$
_{a}^{RL}D^{\mu}_{t}\circ~ _{a}^{\mathbb{T}}D^{-\mu}_{t} 
= ~_{a}^{RL}D^{\mu}_{t}\circ _{a}^{\mathbb{T}}\mathbb{I}^{\mu}_{t} = I.
$$
Again, from Definition~\ref{r25} and Proposition~\ref{p2:b}, we have 
$$
_{a}^{\mathbb{T}}\mathbb{I}^{-\mu}_{t}\circ~ _{a}^{\mathbb{T}}\mathbb{I}^{\mu}_{t} 
=~ ^{RL}_{a}D^{\mu}_{t}\circ~ ^{\mathbb{T}}_{a}\mathbb{I}^{\mu}_{t} = I,
$$
which concludes the proof.
\end{proof}


\section{Conclusions}
\label{sec5}

In this paper, we discussed the nabla fractional derivative on time scales 
in a unified approach by using Gr\"unwald--Letnikov and Riemann--Liouville derivative, 
respectively. Then, we have initiated the generalized definition of nabla derivative 
in fractional order in a pure sense of Riemann--Liouville {and Caputo}. We claim 
that a lot of further work can be completed by using this new idea. 
The aim of formulating the derivative is to solve fractional dynamic equations, 
{stochastic dynamic equations, fuzzy dynamic equations, and one can think 
to extend the concept in a complex dynamic setting}. About applications, 
it has great prospect in mathematical modeling, for example in epidemiology, 
{anomalous diffusion in magnetic resonance imaging \cite{YAW}, 
fractal derivatives modeling \cite{WY}}, and consensus problems 
in time scales on fractional calculus.

\vspace{6pt}

\authorcontributions{All the authors have equal contribution for the preparation of the
article. All authors have read and agreed to the published version of the manuscript.}

\funding{This research was partially funded by the Portuguese 
Foundation for Science and Technology (FCT),
grant number UIDB/04106/2020 (CIDMA).}

\institutionalreview{Not applicable.}

\informedconsent{Not applicable.}

\dataavailability{Not applicable.} 

\acknowledgments{Delfim F. M. Torres is grateful to CIDMA, UIDB/04106/2020.}

\conflictsofinterest{The authors declare that there are no conflicts of interest.} 

\end{paracol}


\reftitle{References}


\end{document}